\documentclass[11pt,bezier,epsf]{amsart}
\usepackage{amsmath,amssymb,amsfonts}
\usepackage{graphicx}
\usepackage{cite}
\usepackage{subfigure}
\newbox\subfigbox
\makeatletter
{\def\caption##1{\gdef\subcapsave{\relax##1}}%
  \let\subcapsave\@empty%
  \setbox\subfigbox\hbox%
  \bgroup}%
{\egroup%
  \subfigure[\subcapsave]{\box\subfigbox}}%
\makeatother

 \newtheorem{thm}{Theorem}[section]

\theoremstyle{definition}
\newtheorem{example}[thm]{Example}
 \newtheorem{defn}[thm]{Definition}
\theoremstyle{remark}
 \newtheorem{rem}[thm]{Remark}
\numberwithin{equation}{section}
\begin{document}

\title[a note on spaces of continuous functions]{a note on spaces of continuous
functions on compact scattered spaces}

\author[F. Naderi]{Fouad Naderi}

\address{Department of Mathematical and Statistical Sciences, University of
Alberta, Edmonton, Alberta, T6G 2G1, Canada.}
\email{naderi@ualberta.ca}


\begin{abstract}
In 1959, Pelczynski and Semadeni proved a theorem in which they gave some
equivalent conditions for a compact Hausdorff space to be scattered.
The purpose of the current note is that to clarify the meaning of the subtle term "conditionally weakly sequentially compact" they used as the basis for the proof of their
theorem. Unfortunately,  the term now is taken over by a similar but subtle concept that may cause a serious problem.
\end{abstract}

\maketitle
%%%%%%%%%%%%%%%%%%%%%%%%%%%%%%%%%%%%%%%%%%%%%%%%%%%%%%%%%%%%%%%%%%%%%%%%%%%%%%%%%%%%%%%%%%%%%

\section{The main result}

 A locally compact Hausdorff topological space $\varOmega$ is {\it scattered} if
 $\varOmega$ does not contain any non-empty perfect subset (i.e., a closed
 non-empty subset $\varPi$ of $\varOmega$ such that each point of $\varPi$ is an
 accumulation point of $\varPi$). Equivalently, any non-empty subset of
 $\varOmega$ contains at least one isolated point. Some authors use the term
 {\it dispersed} instead of scattered. For more details on scattered spaces see
 \cite{Mont}, \cite{Semadeni} and \cite{Semadeni_book}.

\begin{example} \label{infinity}
Consider $\mathbb{N}$ with its usual topology and its one point compactification
${\mathbb{N}}^{*}=\mathbb{N} \cup \{\infty\}$. Both $\mathbb{N}$ and
${\mathbb{N}}^{*}$ are scattered. Also, $\{\frac{1}{n}: n\in \mathbb{N} \}\cup \{0\} $
is another compact scattered space which is not discrete. It can be
shown that a compact metric space is scattered if and only if it is countable
\cite[p.737]{Mont}. A scattered 
space is always totally disconnected. The converse of this is not true as seen by the set $\mathbb{Q}$ of rational numbers.
\end{example}

Consider the following two definitions for the term {\it conditionally weakly
sequentially compact}.

\begin{defn} \label{ES}
Let $\varOmega$ be a compact Hausdorff space. We say that $C(\varOmega)$ has
{\it conditionally weakly sequentially compact property in the sense of\quad {\bf ES}} in the sense that for every
bounded sequence $(x_n)$ of elements of $C(\varOmega)$ there exists a subsequence $(x_{n_{k}})$ and a member $ x_0 \in C(\varOmega)$ such that for every bounded linear functional $\xi$ on
$C(\varOmega)$ the sequence of numbers $\xi(x_{n_{k}})$ converges to $\xi(x_0)$. In other words, if $S$ is a bounded subset of $C(\varOmega)$, then $S$ must be conditionally (=relatively) weakly sequentially compact in the modern language. 
\end{defn}

\begin{defn} \label{PS}
Let $\varOmega$ be a compact Hausdorff space. We say that $C(\varOmega)$ has
{\it conditionally weakly sequentially compact property in the sense of\quad {\bf PS}} if for every
bounded sequence $(x_n)$ of elements of $C(\varOmega)$ there exists a subsequence $(x_{n_{k}})$ such that for every bounded linear functional $\xi$ on
$C(\varOmega)$ the sequence of numbers $\xi(x_{n_{k}})$ is convergent.
\end{defn}
\begin{rem}
The way Definition \ref{PS} seems is now clear, but it was after a correspondence with Professor Semadeni that I could write down it this way. In the Definition \ref{ES}, the space is {\bf weakly sequentially complete} while in the second one we do not need such a strong condition (not to the point $x_0$). Meanwhile, we use Definition \ref{ES} in Eberlin-Smulian theorem to assure weakly compactness of a given set.
\end{rem}

\begin{thm} \label{scatter}
	Let $\varOmega$ be a compact Hausdorff space and $C(\varOmega)$ have conditionally weakly sequentially compact property in the sense of \quad {\bf ES}. Then $\varOmega$ is finite.
\end{thm}
{\bf Proof.} Suppose to the contrary that $\varOmega$ is infinite. Then, the
Banach space of continuous functions $C(\varOmega)$ is infinite dimensional. It
is well-known that the weak closure of the unit sphere of an infinite
dimensional Banach space is the closed unit ball of the space
\cite[p.128]{Conway}. Therefore, if $S$ is the unit sphere of $C(\varOmega)$,
then ${\overline{S}}^{wk}$ is equal to the unit ball $B$ of $C(\varOmega)$.
Since $C(\varOmega)$ has propert \quad {\bf ES}, ${\overline{S}}^{wk}=B$ is weakly
sequentially compact. According to Eberlin-Smulian Theorem \cite[p.163]{Conway},
$B$ must be weakly compact. Therefore by \cite[Theorem 4.2, p.132]{Conway},
$C(\varOmega)$ must be reflexive. By \cite[p.90]{Conway}, $\varOmega$ must be
finite. But, this contradicts our assumption that $\varOmega$ was infinite.
Hence, $\varOmega$ can only be a finite set.$\blacksquare$

\begin{rem}
Suppose $\varOmega$ is a compact Hausdorff space and the main theorem of
\cite[p.214]{Semadeni} holds. Condition (0) and (9) of the latter theorem assert
that $\varOmega$ is scattered if and only if $C(\varOmega)$ has
 conditionally weakly sequentially compact property in the sense of {\bf PS}. But, if one wants to classify scatteredness of $\varOmega$ in the sense of {\bf EP}, (s)he would always end up with a finite set! which is not always the case as Example \ref{scatter} indicates.
\end{rem}

{\bf Acknowledgments.} The author would like to thank Professor Z. Semadeni for his careful comments. He would also thanks Professor W. Zelazko who made this corespondents possible.

%%%%%%%%%%%%%%%%%%%%%%%%%%%%%%%%%%%%%%%%%%%%%%%%%%%%%%%%%%%%%%%%%%%%%%%%%%%%%%%%%%%%%%%%%%%%%%%

%%%%%%%%%%%%%%%%%%%%%%%%%%%%%%%%%%%%%%%%%%%%%%%%%%%%%%%%%%%%%%%%%%%%%%%%%%%%%%%%%%%%%%%%%%%%%%%%%%%%%%%%%%%%%%%%%%%%%%%%%%%%%%
\end{document}